\documentclass[12pt]{amsart}
\usepackage{geometry} 
\usepackage{graphicx}
\usepackage{amsmath}
\usepackage{amsfonts}
\usepackage{mathrsfs}
\usepackage{amssymb}
\usepackage{wrapfig}
\usepackage{bbold}
\usepackage{xcolor, soul}
\usepackage{mdframed}
\usepackage{enumitem}
\usepackage{hyperref}

\newtheorem{thm}{Theorem}
\newtheorem{lem}{Lemma}
\newtheorem{cor}{Corollary}

\newenvironment{pf}{\noindent {\em Proof}.\ \ }{\hspace*{\fill}\rule{.5ex}{1.4ex}\,}

\numberwithin{equation}{section}
\numberwithin{defn}{section}

\newcommand{\real}{\mathbb{R}}
\newcommand{\complex}{\mathbb{C}}

\newcommand{\integer}{\mathbb{Z}}

\newcommand{\one}{\mathbb{1}}

\DeclareMathOperator{\diag}{diag}

\newcommand{\lloc}{L^1_{\mbox{\tiny loc}}(\mathbb{R})}

\newcommand{\acloc}{AC_{\mbox{\tiny loc}}(\mathbb{R})}

\newcommand{\simp}{\mathcal{S}}
\newcommand{\coeff}[1]{\mathbf{#1}}

\title{Integration of an arbitrary linear ODE}

\author{Peter C.~Gibson}
\address{Dept.~of Mathematics \& Statistics\\
York University\\
Toronto, {ON}, M3J1P3, Canada}
\email{pcgibson@yorku.ca}

\date{December 26, 2025} 

\begin{document}

\maketitle

\begin{abstract}
The standard text book theory of ODEs lacks a general method to solve linear equations having variable coefficients, providing instead a collection of special techniques for particular classes of equations.  The present article addresses this shortcoming in the basic theory. 
We introduce the multex integral operator, 
generalizing to several input functions the standard exponential 
primitive operator 
that is inverse to the logarithmic derivative.  
The multex operator serves to integrate in explicit form an arbitrary linear ordinary differential equation.  
\end{abstract}

\vspace*{12pt}
\keywords{\textbf{Key words}: linear ODE; explicit solution }

\subjclass{\textbf{MSC}: 34A05}

\section{Introduction\label{sec-introduction}}

Our goal is to solve explicitly the homogeneous ordinary differential equation
\begin{equation}\label{ode}
y^{(n)}=a_1y^{(n-1)}+a_2y^{(n-2)}+\cdots+a_{n-1}y^\prime+a_ny,
\end{equation}
where coefficients $a_1,\ldots,a_n\in\lloc$ are arbitrary complex valued functions and the unknown function $y$ satisfies $y^{(n-1)}\in\acloc$. The corresponding inhomogeneous equation may then be solved by standard variation of parameters. 
The case $n=2$ of equation (\ref{ode}) was analyzed recently in \cite{Gi:arXiv2025}. The present paper addresses higher order cases, requiring significant
generalization of the methods used for $n=2$.

\subsection{Background\label{sec-background}}

The importance of (\ref{ode}) stems from the wide application of linear ODEs throughout the sciences, both directly as models of physical phenomena, and indirectly, for example in the analysis of nonlinear PDEs by scattering theory, such as applies to the Gelfand-Dikii hierarchy \cite{BeDeTo:1988}. 
 Conventional wisdom as expressed in the text book literature, including the excellent up-to-date reference \cite[\S3.5]{Te:2012}, asserts that there is no general way of solving $n$th order equations. This viewpoint seems reasonable in light of the antiquity of the subject matter (e.g., variation of parameters dates to Lagrange \cite{La:1869}) and its standard presentation \cite{In:1944,BoDi:1965,encyclopedia_mathematics_linear_ode}, but the assertion is not entirely accurate. 
 
Given a locally integrable matrix valued function $\mathbf{m}:\real\rightarrow\complex^{n\times n}$, the first order initial value problem in $n\times n$ matrices, 
\begin{equation}\label{system}
M^\prime=\coeff{m}M,\qquad M(0)=I
\end{equation}
has an explicit locally absolutely continuous solution for $x>0$ of the form
\begin{multline}\label{system-solution}
M(x)=I+\int\limits_{0<s_1<x}\coeff{m}(s_1)\,ds_1+\\
\int\limits_{0<s_1<s_2<x}\coeff{m}(s_2)\coeff{m}(s_1)\,ds_1ds_2+\int\limits_{0<s_1<s_2<s_3<x}\coeff{m}(s_3)\coeff{m}(s_2)\coeff{m}(s_1)\,ds_1ds_2ds_3+\cdots
\end{multline}
See Appendix~\ref{sec-convergence}. Such formulations appear in the work of Peano \cite{Pe:1888} and, more fully fledged, in that of Volterra and Hostinsky \cite[Ch.VI,\S2]{HoVo:1938}. In physics related literature and elsewhere, the solution $M$ to (\ref{system}) is known as the ordered exponential (or product integral) of the coefficient matrix $\coeff{m}$ \cite{DoFr:1977,GiJo:1990,Sl:2007}.  If $\coeff{m}$ belongs to a particular Lie algebra, then the solution $M$ belongs to the corresponding Lie group. 
Note that in the scalar case $n=1$, the right-hand side of (\ref{system-solution}) reproduces the usual exponential primitive $e^{\int_0^x\coeff{m}}$. When $\mathbf{m}$ is a Hamiltonian operator, the right-hand side of (\ref{system-solution}) is known as a Dyson series (see \cite{Dy:1949} and \cite[\S5.7.1]{SaNa:2021}), a term which will also be used here for the matrix case. 

Since the $n$th order equation (\ref{ode}) can be represented as a first-order system, the Dyson series gives an explicit formulation of the solution.  However, in the case of a single equation, the Dyson series yields an $n\times n$ matrix, whereas all that is needed is an $n$-element basis for the general solution.  This $n$-fold redundancy raises a basic question: is there is a purely scalar formulation of the solution to (\ref{ode}) that does not entail prior determination of $n^2$ functions as in (\ref{system-solution})? The present paper is devoted precisely to this question.  We demonstrate an affirmative answer, with the resulting formulation valid locally about any given point. Without loss of generality, the point in question may be taken to be $x=0$, since the general form of (\ref{ode}) is invariant under translation. 
As an illustration of our main result, consider the third order initial value problem
\begin{equation}\label{third-order-example}
y^{(3)}=\alpha y^{\prime}+\beta y,\qquad y(0)=1,\quad y^{\prime}(0)=y^{\prime\prime}(0)=0 
\end{equation}
where $\alpha,\beta\in\lloc$ are arbitrary.  For $x>0$, set 
\begin{equation}\label{gamma}
\gamma(x)=1+\negthickspace\negthickspace\int\limits_{0<s_1<s_2<x}\negthickspace\negthickspace\alpha(s_1)\,ds_1ds_2+\negthickspace\negthickspace\negthickspace\negthickspace
\int\limits_{0<s_1<s_2<s_3<s_4<x}\negthickspace\negthickspace\negthickspace\negthickspace\alpha(s_1)\alpha(s_3)\,ds_1ds_2ds_3ds_4+\cdots
\end{equation}
The series (\ref{gamma}) converges absolutely, uniformly on compact sets, since $\alpha\in\lloc$. Moreover, $\gamma^\prime(x)$ is locally absolutely continuous for all $x>0$, and $\gamma(0+)=1$.  Indeed, it is straightforward to verify that $\gamma$ solves the initial value problem
\begin{equation}\label{gamma-ivp}
\gamma^{\prime\prime}=\alpha\gamma,\qquad\gamma(0)=1,\quad\gamma^\prime(0)=0. 
\end{equation}
Denote
\begin{equation}\label{Phi}
\Phi(s_1,s_2,s_3)=\frac{\beta(s_1)\gamma(s_1)\gamma(s_3)}{\gamma(s_2)^2}.
\end{equation}
Then for all sufficiently small $x>0$, the solution to (\ref{third-order-example}) is given by the absolutely convergent integral series 
\begin{multline}\label{third-order-solution}
y(x)=1+\negthickspace\negthickspace\int\limits_{0<s_1<s_2<s_3<x}\negthickspace\negthickspace\Phi(s_1,s_2,s_3)\,ds_1ds_2ds_3\\
+\negthickspace\negthickspace\int\limits_{0<s_1<\cdots<s_6<x}\negthickspace\negthickspace\Phi(s_1,s_2,s_3)\Phi(s_4,s_5,s_6)\,ds_1\cdots ds_6+\cdots
\end{multline}
A similar formula holds for all $x<0$ sufficiently close to 0. 

Note that (\ref{third-order-solution}) differs from a Dyson series in that the integrands are scalar functions, not products of matrices or non-commuting operators. Furthermore, the domain of each integral is a simplex whose dimension is congruent to 0 modulo 3, whereas a Dyson series involves every dimension. 
Remarkably, \emph{any} equation of the form (\ref{ode}), with arbitrary initial conditions, can be solved explicitly in terms of series akin to (\ref{third-order-solution}).  Our main result is a precise formulation of this fact. 

A feature of the solution (\ref{third-order-solution}) is that it involves the solution $\gamma$ to an auxiliary equation (\ref{gamma-ivp}) of lower order than the original (\ref{third-order-example}).  This is a general phenomenon.  The solution to any equation of the form (\ref{ode}) can be formulated in terms of what we call auxiliary functions, i.e., solutions to linear initial value problems of strictly lower order.  

\subsection{Plan of the paper\label{sec-plan}}

The present paper is essentially self contained. For ease of reference, a brief justification of the Dyson formula (\ref{system-solution}) in the case where the coefficient matrix $\coeff{m}$ is entrywise locally integrable is included in Appendix~\ref{sec-convergence}. 

In view of the scope of (\ref{ode}), our main result is strikingly concise;
the explicit general solution to (\ref{ode}) can be formulated in just a few lines. Two preliminary constructions make this possible. The first is that of the multex operator, given in \S\ref{sec-multex-operator}, along with its particular summands called trig operators, described in \S\ref{sec-trig-operators}. The second construction, of auxiliary functions determined by the coefficients of (\ref{ode}), is detailed in \S\ref{sec-auxiliary}. 

Our main results are Theorem~\ref{thm-main} and Corollary~\ref{cor-examples}, stated in \S\ref{sec-main-result}. 
Corollary~\ref{cor-examples} applies the main theorem iteratively to equations of order $n=2,3$ and $4$, rendering explicit the various auxiliary functions involved.  The passage from Theorem~\ref{thm-main} to Corollary~\ref{cor-examples} is purely mechanical, and in principle it is completely straightforward. Nevertheless, the first steps of this passage are detailed in Appendix \ref{sec-derivation} to clarify 
the intended interpretation of notation etc.  

Theorem~\ref{thm-main} is very widely applicable, and no attempt has been made to describe applications exhaustively.  Two illustrative examples, the Schr\"odinger and Sommerfeld-Orr equations, are included in \S\ref{sec-examples}. 
Both are of independent interest; analysis of the former led to a precursor of Theorem~\ref{thm-main} in \cite{Gi:Proc2026}.

Theorem~\ref{thm-main} is proved in \S\ref{sec-proof}. The proof is broken down into two parts, in \S\ref{sec-proof-part-one} and \S\ref{proof-part-two}, according to initial values. 

The final section of the paper, \S\ref{sec-remarks}, includes a short discussion of the results. 

\section{Preliminary constructions\label{sec-preliminary}}

\subsection{The multex operator\label{sec-multex-operator}}

For fixed $n\geq1$, let $\nu:\integer_{\geq0}\rightarrow \{1,\ldots,n\}$ map each $j\geq 0$ to the unique element of $\{1,\ldots, n\}$ congruent to $j$ modulo $n$. Define the $n$-ary multivariate exponential primitive operator $$E:\bigl(\lloc)^n\rightarrow\acloc$$ by the formula:
\begin{equation}\label{multex}
\begin{split}
E_{f_1,\ldots,f_n}&(x)=\\
&\left\{
\begin{array}{cc}
1+\sum\limits_{j=1}^\infty\thickspace\displaystyle\int\limits_{0<s_1<\cdots<s_j<x}f_{\nu(1)}(s_1)f_{\nu(2)}(s_2)\cdots f_{\nu(j)}(s_j)\,ds_1\cdots ds_j&\mbox{ if }x\geq 0\\
1+\sum\limits_{j=1}^\infty\rule{0pt}{21pt}(-1)^j\negthickspace\negthickspace\negthickspace\displaystyle\int\limits_{x<s_j<\cdots<s_1<0}f_{\nu(1)}(s_1)f_{\nu(2)}(s_2)\cdots f_{\nu(j)}(s_j)\,ds_1\cdots ds_j&\mbox{ if }x\leq 0
\end{array}
\right..
\end{split}
\end{equation}
Observe that if $f_1=f_2=\cdots=f_n$, then $$E_{f_1,\ldots,f_n}(x)=e^{\int_{0}^xf_1}$$ whereby (\ref{multex}) generalizes the usual exponential primitive, inverse to the logarithmic derivative.  Standard arguments show the integral (\ref{multex}) converges absolutely, uniformly on compact sets, and $E_{f_1,\ldots,f_n}(x)$ is locally absolutely continuous; for details see \cite[Lem.~1]{Gi:arXiv2025}.  To avoid repetition of the cumbersome phrase ``$n$-ary multivariate exponential primitive," the operator $(f_1,\ldots,f_n)\mapsto E_{f_1,\ldots,f_n}$ will henceforth be referred to as the \emph{multex} operator, the arity $n$ being clear from context.   

In a neighbourhood of $x=0$, the general solution to (\ref{ode}) can be formulated in terms of the coefficients $a_j$ using various summands of the multex operator (\ref{multex}).

\subsection{Trig operators\label{sec-trig-operators}}

The multex operator is comprised of simplicial integrals 
\begin{equation}\label{simplicial}
\simp^j_{f_1,\ldots,f_n}(x)=
\left\{
\begin{array}{cc}
\displaystyle\int\limits_{0<s_1<\cdots<s_j<x}f_{\nu(1)}(s_1)f_{\nu(2)}(s_2)\cdots f_{\nu(j)}(s_j)\,ds_1\cdots ds_j&\mbox{ if }x\geq 0\\[20pt]
(-1)^j\negthickspace\negthickspace\negthickspace\displaystyle\int\limits_{x<s_j<\cdots<s_1<0}f_{\nu(1)}(s_1)f_{\nu(2)}(s_2)\cdots f_{\nu(j)}(s_j)\,ds_1\cdots ds_j&\mbox{ if }x\leq 0
\end{array}
\right.
\end{equation}
where $j\geq1$ corresponds to the dimension of the underlying simplex 
\[
0<s_1<\cdots<s_j<x\quad\mbox{ or }\quad x<s_j<\cdots<s_1<0.
\]
Denote $\simp^0_{f_1,\ldots,f_n}(x)=1$. For each $j$ in the range $1\leq j\leq n$, define
\begin{equation}\label{T-0}
T_{f_1,\ldots,f_n;j}=\sum_{\stackrel{k\geq 0}{\nu(k)=j}}\simp^k_{f_1,\ldots,f_n}.
\end{equation}

Alternatively, 
set $\varepsilon_{j,k}^n=-1$ if $k$ is congruent to $j$ or $j+1$ modulo $n$, and $\varepsilon_{j,k}^n=1$ otherwise.  
Then, as is easily verified, one may equivalently define
\begin{equation}\label{T}
T_{f_1,\ldots,f_n;j}=\left\{
\begin{array}{cl}
\textstyle\frac{1}{2}\bigl(E_{f_1,\ldots,f_n}-E_{\varepsilon_{j,1}^nf_1,\ldots,\varepsilon_{j,n}^nf_n}\bigr)&\mbox{ if }1\leq j\leq n-1\\[5pt]
\textstyle\frac{1}{2}\bigl(E_{f_1,\ldots,f_n}+E_{\varepsilon_{j,1}^nf_1,\ldots,\varepsilon_{j,n}^nf_n}\bigr)&\mbox{ if }j=n
\end{array}\right..
\end{equation}

It follows from (\ref{multex}) and  (\ref{T-0}) (or from (\ref{multex}) and (\ref{T})), that 
\begin{equation}\label{decomposition}
E_{f_1,\ldots,f_n}(x)=\sum_{j=1}^nT_{f_1,\ldots,f_n;j}(x)
\end{equation}
and furthermore that
\begin{equation}\label{zero-values}
T_{f_1,\ldots,f_n;j}(0)=\delta_{j,n},
\end{equation}
where $\delta_{j,n}$ is the Kronecker delta.  The decomposition (\ref{decomposition}) of the multex operator as a sum of component operators is reminiscent of the elementary decomposition
\[
e^x=\cosh x+\sinh x,
\]
just as the definition (\ref{T}) is analogous to that of $\cosh$ and $\sinh$ in terms of $\exp$. 
Thus viewing the multex operator $(f_1,\ldots,f_n)\mapsto E_{f_1,\ldots,f_n}$ as a sort of exponential, the corresponding operators $(f_1,\ldots,f_n)\mapsto T_{f_1,\ldots,f_n;j}$ as defined by (\ref{T-0}) or (\ref{T}) may  be viewed as generalized trigonometric operators.  For brevity we shall refer to them simply as \emph{trig} operators. 
The trig operators satisfy the following key differential relations. 
\begin{lem}\label{lem-trig}
\begin{equation}\label{trig-derivatives}
\frac{d}{dx}T_{f_1,\ldots,f_n;j}(x)=\left\{
\begin{array}{cl}
f_j(x)T_{f_1,\ldots,f_n;j-1}(x)&\mbox{ if }\quad 2\leq j\leq n\\[5pt]
f_j(x)T_{f_1,\ldots,f_n;n}(x)&\mbox{ if }\quad j=1
\end{array}\right.
\end{equation}
\end{lem}
\begin{pf} The formula (\ref{simplicial}) implies that for any $j\geq1$, 
\begin{equation}\label{simp-derivative}
\frac{d}{dx}\simp^j_{f_1,\ldots,f_n}(x)=f_{\nu(j)}(x)\simp^{j-1}_{f_1,\ldots,f_n}(x). 
\end{equation}
Equation (\ref{trig-derivatives}) then follows by (\ref{T-0}).
\end{pf}

\subsection{Auxiliary functions\label{sec-auxiliary}}

Let $n\geq 2$, fix $a_1,\ldots,a_n\in\lloc$, and set 
\begin{equation}\nonumber
a=(-1,a_1,\ldots,a_n)^t.
\end{equation}
Define associated \emph{auxiliary functions} $\varphi_1,\ldots,\varphi_n$ in a neighbourhood of $x=0$, starting with $\varphi_n$ and iterating downwards, as follows.  Denote $D=\frac{d}{dx}$, and set 
\begin{equation}\label{D-matrix}
\mathscr{D}=
\begin{pmatrix}0&1&D&D^2&\cdots&D^{n-1}\\
0&0&1&D&\cdots&D^{n-2}\\
0&0&0&1&\ddots&\vdots\\
\vdots&&&\ddots&\ddots&D\\
&&&&&1\\
0&&&\cdots&&0
\end{pmatrix}.
\end{equation}
Thus $\mathscr{D}$ is an $(n+1)\times(n+1)$ Toeplitz matrix, with first row $(0,1,D,D^2,\cdots,D^{n-1})$. For each $1\leq k\leq n$ denote $\alpha^k=(\alpha^k_1,\ldots,\alpha^k_{n+1})^t$. Set 
\begin{equation}\nonumber
\alpha^n=(0,\cdots,0,1)^t\in\complex^{n+1}.
\end{equation}
For each $2\leq k\leq n$, given $\alpha^{k}$, let $\varphi_k$ denote the unique local solution to the \emph{auxiliary equation} in unknown scalar function $u$, 
\begin{equation}\label{k-ode}
a^t\mathscr{D}\bigl(u\alpha^{k}\bigr)=0,
\end{equation}
such that $u=\varphi_k$ satisfies initial conditions
\begin{equation}\label{k-initial}
u^{(j-1)}(0)=\delta_{j,1}\qquad(1\leq j\leq k-1).
\end{equation}
Set 
\begin{equation}\label{alpha-k}
\alpha^{k-1}=\mathscr{D}\bigl(\varphi_k\alpha^{k}\bigr).
\end{equation}

Observe by (\ref{k-ode}) and (\ref{alpha-k}) that the equation determining $\varphi_n$ has order $n-1$, and that $\alpha^{n-1}_{n+1}=0$.  For $2\leq k\leq n-1$, if $\alpha^k_j=0$ for each $j$ in the range $k+1\leq j\leq n+1$, then it follows from (\ref{k-ode}) that the auxiliary equation defining $\varphi_k$ has order at most $k-1$.  And (\ref{alpha-k}) implies that $\alpha^{k-1}_j=0$ for every $j$ in the range $k\leq j\leq n+1$.  Therefore by downward induction, the auxiliary equation (\ref{k-ode}) has order at most $k-1$ for all $2\leq k\leq n$, and $\alpha^k_j=0$ for every $k+1\leq j\leq n+1$.  

Starting with the constant vector $\alpha^{n}=(0,\ldots,0,1)^t\in\complex^{n+1}$, one may thus iteratively determine $\varphi_{n},\varphi_{n-1},\ldots,\varphi_2$ in a neighbourhood of $x=0$ by solving equations of progressively lower order. Lastly, define $\varphi_1$ by the equation 
\begin{equation}\label{phi-1-2}
\varphi_1\varphi_2\cdots\varphi_n=a_n. 
\end{equation}

\section{Main results\label{sec-main-result}}

Denote by $\pi_n$ the right shift permutation, for which $\pi_n(j)$ is the unique element of $\{1,\ldots,n\}$ congruent to $j-1$ modulo $n$. Let $\pi_n^k$ denote its $k$-fold composition, with $\pi_n^0$ being the identity. 

\begin{thm}\label{thm-main}
Let $n\geq 2$, and fix $a_1,\ldots,a_n\in\lloc$. Let $\varphi_1,\ldots,\varphi_n$ be the associated auxiliary functions. For each $1\leq k\leq n$, denote by $\psi_{n,k}$ the unique solution to 
\begin{equation}\label{ode-n}
y^{(n)}=a_1y^{(n-1)}+a_2y^{(n-2)}+\cdots+a_{n-1}y^\prime+a_ny
\end{equation}
such that $y=\psi_{n,k}$ satisfies initial conditions
\begin{equation}\label{ic-j}
y^{(j-1)}(0)=\delta_{j,k}\qquad(1\leq j\leq n).
\end{equation}
Then for each $1\leq k\leq n$, there is an open neighbourhood of $x=0$ in which $\psi_{n,k}(x)$ is given explicitly in terms of the auxiliary functions by the formula
\begin{equation}\label{psi-n}
\psi_{n,k}=T_{\varphi_{\eta(1)},\ldots,\varphi_{\eta(n)};\eta^{-1}(n)}\quad\mbox{ where }\quad\eta=\pi^{k-1}_n.
\end{equation}
\end{thm}
The functions $\psi_{n,k}$ are defined such that the general solution to (\ref{ode-n}) has the form
\begin{equation}\label{general-solution}
y=\sum_{k=1}^ny^{(k-1)}(0)\psi_{n,k}.
\end{equation}
Thus in a neighbourhood of $x=0$ the general solution to (\ref{ode-n}) is expressed in terms of solutions to lower order auxiliary equations by means of trig operators.  An explicit solution in terms of the coefficients $a_j$ is obtained by recursive application of the theorem to the auxiliary equations, as the following result illustrates. 

\begin{cor}\label{cor-examples} Let $P:\lloc\rightarrow\acloc$ denote the primitive operator $$Pf(x)=\int_{0}^xf.$$
Fix $a_1,a_2,a_3,a_4\in\lloc$.\\[0pt] 
\begin{enumerate}[label={\textbf{(\Roman*)}},itemindent=0em]
\item\label{item-second-order}
Set 
\begin{equation}\nonumber
C=T_{a_2e^{-Pa_1},e^{Pa_1};2}\mbox{ and }S=T_{e^{Pa_1},a_2e^{-Pa_1};1}.
\end{equation}
The global solution to the second order equation 
\begin{equation}\nonumber
y^{\prime\prime}=a_1y^\prime+a_2y
\end{equation}
is given by the formula
\begin{equation}\nonumber
y=y(0)C+y^{\prime}(0)S.
\end{equation}
\item\label{item-third-order}
Set
\begin{equation}\label{phi-3-0}
\varphi_3=C=T_{a_2e^{-Pa_1},e^{Pa_1};2},\quad\varphi_2=e^{Pa_1}C^{-2},\quad\mbox{ and }\quad\varphi_1=a_3e^{-Pa_1}C.
\end{equation}
In a neighbourhood of $x=0$, the general solution to the third order equation
\begin{equation}\nonumber
y^{(3)}=a_1y^{\prime\prime}+a_2y^\prime+a_3y
\end{equation}
is given by the formula
\begin{equation}\nonumber
y=y(0)T_{\varphi_1,\varphi_2,\varphi_3;3}+y^\prime(0)T_{\varphi_3,\varphi_1,\varphi_2;1}+y^{\prime\prime}(0)T_{\varphi_2,\varphi_3,\varphi_1;2}.
\end{equation}
\item\label{item-fourth-order}
Referring to (\ref{phi-3-0}), set 
\begin{equation}\nonumber
\begin{split}
\psi_4&=T_{\varphi_1,\varphi_2,\varphi_3;3}\qquad\quad \psi_3=T_{e^{-Pa_1}\psi_4^2(a_2\psi_4+2a_1\psi_4^\prime-3\psi_4^{\prime\prime}),\,e^{Pa_1}\psi_4^{-3};2}\\
\psi_2&=e^{Pa_1}\psi_3^{-2}\psi_4^{-3}\qquad \psi_1=a_4e^{-Pa_1}\psi_3\psi_4^2.
\end{split}
\end{equation}
In a neighbourhood of $x=0$, the general solution to the fourth order equation
\begin{equation}\nonumber
y^{(4)}=a_1y^{(3)}+a_2y^{\prime\prime}+a_3y^\prime+a_4y
\end{equation}
is given by the formula
\begin{equation}\nonumber
\begin{split}
y=y(0)T_{\psi_1,\psi_2,\psi_3,\psi_4;4}+&y^\prime(0)T_{\psi_4,\psi_1,\psi_2,\psi_3;1}+\\ &y^{\prime\prime}(0)T_{\psi_3,\psi_4,\psi_1,\psi_2;2}+y^{(3)}(0)T_{\psi_2,\psi_3,\psi_4,\psi_1;3}.
\end{split}
\end{equation}
\end{enumerate}
\end{cor}

 Corollary~\ref{cor-examples} presents explicit solutions to the equations most prevalent in applications. 
However, the iterative procedure described in Theorem~\ref{thm-main} may be continued indefinitely to arbitrary order. 

The conclusion in part~\ref{item-second-order} that the solution is global follows from global existence of the auxiliary functions in the case $n=2$.  See \cite{Gi:arXiv2025} for further details.  For $n>2$ the auxiliary functions may have singularities, as evident in the formula (\ref{Phi}) occurring in the introductory example of \S\ref{sec-background}. 

\section{Examples: the Schr\"odinger and Orr-Sommerfeld equations\label{sec-examples}}

\subsection{The 1D Schr\"odinger equation in impedance form\label{sec-1d}}
Given a positive, real-valued function $\zeta$ and a constant $\omega$, the equation 
\begin{equation}\label{helmholtz}
(\zeta u^\prime)^\prime+\omega^2\zeta u=0
\end{equation}
is known as the Schr\"odinger equation in impedance form, or the Helmholtz equation.  Simple transformations relate it both to the classical Schr\"odinger equation and the variable coefficient wave equation.  Provided $\zeta$ is locally absolutely continuous,  (\ref{helmholtz}) can be written in the form
\begin{equation}\label{helmholtz-2}
u^{\prime\prime}=-\textstyle\frac{\zeta^\prime}{\zeta}u^\prime-\omega^2u
\end{equation}
to which Corollary~\ref{cor-examples} applies with $a_1=-\zeta^\prime/\zeta$ and $a_2=-\omega^2$.  Noting that 
\begin{equation}\nonumber
e^{-Pa_1}=\zeta/\zeta(0)\quad\mbox{ and }\quad e^{Pa_1}=\zeta(0)/\zeta
\end{equation}
one thereby obtains the general solution 
\begin{equation}\label{helmholtz-solution}
u=u(0)C+u^\prime(0)S
\end{equation}
where
\begin{equation}\label{C-helmholtz}
\begin{split}
&C(x)=\\
&\left\{\begin{array}{cc}
1+\displaystyle\sum_{n=1}^\infty(i\omega)^{2n}\negthickspace\negthickspace\int\limits_{0<s_1<\cdots<s_{2n}<x}\negthickspace\negthickspace\frac{\zeta(s_1)\zeta(s_3)\cdots\zeta(s_{2n-1})}{\zeta(s_2)\zeta(s_4)\cdots\zeta(s_{2n})}\,ds_1\cdots ds_{2n}&\mbox{ if }x\geq0\\[25pt]
1+\displaystyle\sum_{n=1}^\infty(i\omega)^{2n}\negthickspace\negthickspace\int\limits_{x<s_{2n}<\cdots<s_1<0}\negthickspace\negthickspace\frac{\zeta(s_1)\zeta(s_3)\cdots\zeta(s_{2n-1})}{\zeta(s_2)\zeta(s_4)\cdots\zeta(s_{2n})}\,ds_1\cdots ds_{2n}&\mbox{ if }x\leq0
\end{array}\right.
\end{split}
\end{equation}
and
\begin{equation}\label{S-hemholtz}
\begin{split}
S(x)&=\\
&\left\{\begin{array}{cc}
\zeta(0)\displaystyle\sum_{n=1}^\infty(i\omega)^{2n-2}\negthickspace\negthickspace\negthickspace\negthickspace\int\limits_{0<s_1<\cdots<s_{2n-1}<x}\negthickspace\negthickspace\frac{\zeta(s_2)\zeta(s_4)\cdots\zeta(s_{2n-2})}{\zeta(s_1)\zeta(s_3)\cdots\zeta(s_{2n-1})}\,ds_1\cdots ds_{2n-1}&\mbox{ if }x\geq0\\[25pt]
-\zeta(0)\displaystyle\sum_{n=1}^\infty(i\omega)^{2n-2}\negthickspace\negthickspace\negthickspace\negthickspace\int\limits_{x<s_{2n-1}<\cdots<s_1<0}\negthickspace\negthickspace\frac{\zeta(s_2)\zeta(s_4)\cdots\zeta(s_{2n-2})}{\zeta(s_1)\zeta(s_3)\cdots\zeta(s_{2n-1})}\,ds_1\cdots ds_{2n-1}&\mbox{ if }x\leq0
\end{array}\right..
\end{split}
\end{equation}
The above formulas reduce easily to the elementary functions $C(x)=\cos(\omega x)$ and $S(x)=\frac{1}{\omega}\sin(\omega x)$ if $\zeta$ is constant. 
The solution to (\ref{helmholtz}) was recently given an explicit formulation in \cite{Gi:SIMA2024}. Whereas the latter involves a lengthy derivation, Corollary~\ref{cor-examples} gives a much more immediate result, which is also simpler.  

\subsection{The Orr-Sommerfeld equation\label{sec-orr-sommerfeld}}

Analysis of hydrodynamic stability of a two-dimensional flow involves the Orr-Sommerfeld equation, which has the generic form
\begin{equation}\label{orr}
y^{(4)}=a_2y^{\prime\prime}+a_4y.
\end{equation}
There is no known analytic solution to (\ref{orr}); see \cite{He:1987,LaTr:1997,GrNg:2019,PiMo:2025}.  Application of Corollary~\ref{cor-examples} gives one, as follows. 

Denote by $\one$ the function with constant value 1, and, referring to the coefficients $a_2,a_4$ in (\ref{orr}),  set 
\begin{equation}\label{C-orr}
C=T_{a_2,\one;2}. 
\end{equation}
The solution to (\ref{orr}) near $x=0$ as prescribed by Corollary~\ref{cor-examples} then works out to be 
\begin{multline}\label{orr-solution}
y=y(0)T_{a_4C,C^{-2},C,\one;4}+y^\prime(0)T_{\one,a_4C,C^{-2},C;1}\\
+y^{\prime\prime}(0)T_{C,\one,a_4C^{-2};2}+y^{(3)}(0)T_{C^{-2},C,\one,a_4C;3}.
\end{multline}

\section{Proof of the main results\label{sec-proof}} 

\subsection{Proof of Theorem~\ref{thm-main}, part 1\label{sec-proof-part-one}}

The goal in this section is to prove that $\psi_{n,1}$, as formulated in the statement of Theorem~\ref{thm-main}, is indeed a solution to equation (\ref{ode-n}) consistent with the given initial conditions.  That the remaining $\psi_{n,k}$ are also solutions will be proved in the next section. 

Fix $a_1,\ldots, a_n$, and set $a_0=-1$.  Denote by $\varphi_1,\ldots,\varphi_n$ the corresponding auxiliary functions as defined in \S\ref{sec-auxiliary}. Denote 
\begin{equation}\nonumber
\tau_j=T_{\varphi_1,\ldots,\varphi_n;j}
\end{equation}
so that by Lemma~\ref{lem-trig}, $\tau_1^\prime=\varphi_1\tau_n$ and 
\begin{equation}\nonumber
\tau_j^\prime=\varphi_j\tau_{j-1}\qquad(2\leq j\leq n). 
\end{equation}
Denote by $L$ the $n\times n$ cyclic left shift permutation matrix
\begin{equation}\nonumber
L=\begin{pmatrix}0&1&0&0&0&\cdots&0\\ 
0&0&1&0&0&\cdots&0\\
0&0&0&1&0&\cdots&0\\
\vdots&&&\ddots&\ddots&&\vdots\\
0&&&&&&1\\
1&0&0&&\cdots&&0
\end{pmatrix}.
\end{equation}
Set 
\begin{equation}\nonumber
D=\frac{d}{dx}\quad\mbox{ and }\quad \Phi=L\diag(\varphi_1,\ldots,\varphi_n)=\begin{pmatrix}0&\varphi_2&0&0&0&\cdots&0\\ 
0&0&\varphi_3&0&0&\cdots&0\\
0&0&0&\varphi_4&0&\cdots&0\\
\vdots&&&\ddots&\ddots&&\vdots\\
0&&&&&&\varphi_n\\
\varphi_1&0&0&&\cdots&&0
\end{pmatrix}.
\end{equation}
Observe by Lemma~\ref{lem-trig} that if 
\begin{equation}\nonumber
f=\sum_{j=1}^nc_j\tau_j
\end{equation}
for some coefficient functions $c_j\in\acloc$, then 
\begin{equation}\label{iteration}
f^\prime=\sum_{j=1}^n\gamma_j\tau_j\quad\mbox{ where }\quad \begin{pmatrix}\gamma_1\\ \vdots\\ \gamma_n\end{pmatrix}=(D+\Phi)\begin{pmatrix}c_1\\ \vdots\\ c_n\end{pmatrix}. 
\end{equation}
Denote by $\mathbf{e}^n$ the $n$th standard basis vector $\mathbf{e}^n=(0,\ldots,0,1)^t\in\complex^n$. Iteration of (\ref{iteration}) $j$ times starting with $f=\tau_n$ yields
\begin{equation}\label{tau-derivatives}
\tau_n^{(j)}=(\tau_1,\ldots,\tau_n)(D+\Phi)^{j}\mathbf{e}^n.
\end{equation}
Therefore
\begin{equation}\label{tau-formula}
\sum_{j=0}^na_j\tau_n^{(n-j)}=(\tau_1,\ldots,\tau_n)v
\end{equation}
where
\begin{equation}\label{v}
v=\begin{pmatrix}v_1\\ \vdots\\ v_n\end{pmatrix}=\sum_{j=0}^na_j(D+\Phi)^{n-j}\mathbf{e}^n.
\end{equation}
It follows from (\ref{tau-formula}) that if $v=0$, then necessarily $y=\tau_n$ satisfies equation (\ref{ode-n}). 

The next step is to show that construction of the auxiliary functions as in \S\ref{sec-auxiliary} forces $v=0$.  Set 
\begin{equation}\label{a}
a=(a_0,a_1,\ldots,a_n)^t
\end{equation}
and write equation (\ref{v}) in the form 
\begin{equation}\label{A}
v=Aa\quad\mbox{ where }\quad A=\left((D+\Phi)^n\mathbf{e}^n,(D+\Phi)^{n-1}\mathbf{e}^n,\ldots,(D+\Phi)\mathbf{e}^n,\mathbf{e}^n\right).
\end{equation}
The following observations concerning monomials in $D$ and $\Phi$ facilitate analysis of the $n\times(n+1)$ coefficient matrix $A$ in (\ref{A}).   

Fix $k$ in the range $1\leq k\leq n$, and let $w=w_1w_2\cdots w_k$, where $w_j\in\{D,\Phi\}$ for every $1\leq j\leq k$.  Call $k$ the \emph{total degree} of $w$. 
Let  $\rho_\Phi(w)$ denote the number of indices $j$ such that $w_j=\Phi$, referred to as the \emph{degree of $w$ with respect to $\Phi$}, and let $\rho_D(w)$ denote the number of indices $j$ such that $w_j=D$, so that 
\begin{equation}\nonumber
k=\rho_\Phi(w)+\rho_D(w).
\end{equation}
The action of $w=w_1\cdots w_k$ on a function $F$ with values in $\complex^n$ is to be understood as successive composition from right to left, 
\begin{equation}\nonumber
wF=\bigl(w_1\bigl(w_2\bigl(\cdots\bigl(w_kF\bigr)\cdots\bigr)\bigr)\bigr).
\end{equation}
Observe that $w\mathbf{e}^n\neq 0$ only if $w_k=\Phi$. Moreover, for $1\leq j\leq n-1$,
\begin{equation}\label{phi-degree}
w\mathbf{e}^n\propto\mathbf{e}^j\quad\Leftrightarrow\quad \rho_\Phi(w)=n-j. 
\end{equation}
Note also that 
\begin{equation}\label{phi-n}
\Phi^n\mathbf{e}^n=\varphi_1\varphi_2\cdots\varphi_n\mathbf{e}^n.
\end{equation}
Let $A_{jk}$ denote the $(j,k)$-entry of $A$. Observe by (\ref{A}) that for $1\leq k\leq n$,  $A_{jk}$ is the sum of all terms of the form $\bigl(\mathbf{e}^j\bigr)^tw\mathbf{e}^n$, where $w=w_1\cdots w_{n+1-k}$ is a monomial in $D$ and $\Phi$ such that $w_{n+1-k}=\Phi$, and $\rho_{\Phi}(w)=n-j$.  The set of such terms is non-empty only if $n+1-k\geq n- j$, or equivalently, if $k\leq j+1$; otherwise $A_{jk}=0$.  In the case where $k=j+1$, $A_{jk}$ consists of a single term 
\begin{equation}\nonumber
A_{j\, j+1}=\bigl(\mathbf{e}^j\bigr)^t\Phi^{n-j}\mathbf{e}^n.
\end{equation}

It follows from these observations that the coefficient matrix $A$ in (\ref{A}) has the general $n\times(n+1)$ form
\begin{equation}\label{A-form}
A=\begin{pmatrix}
\ast&\varphi_2\varphi_3\cdots\varphi_n&0&0&0&\cdots&0\\
\ast&\ast&\varphi_3\cdots\varphi_n&0&0&\cdots&0\\
\ast&\ast&\ast&\varphi_4\cdots\varphi_n&0&\cdots&0\\
\vdots&\vdots&\ddots&\ddots&\ddots&\ddots&\vdots\\
D^{n-1}\varphi_n&D^{n-2}\varphi_n&D^{n-3}\varphi_n&\cdots&D\varphi_n&\varphi_n&0\\
\varphi_1\cdots\varphi_n&0&0&0&\cdots&0&1
\end{pmatrix}.
\end{equation}
This implies in particular that
\begin{equation}\label{n-minus-one}
v_n=-\varphi_1\cdots\varphi_n+a_n\quad\mbox{ and }\quad v_{n-1}=\sum_{j=0}^{n-1}a_jD^{n-j-1}\varphi_n.
\end{equation}
The foregoing observations imply furthermore that each row of $A$ apart from the last may be constructed in a simple way from the row below, as follows. 
Set
\begin{equation}\label{alpha-k-2}
\alpha^k=\bigl(A_{k 1},\ldots,A_{k\, n+1}\bigr)^t. 
\end{equation}
Then for every $2\leq k\leq n$, 
\begin{equation}\label{iteration-2}
\alpha^{k-1}=\mathscr{D}\bigl(\varphi_k\alpha^k\bigr)\quad\mbox{ where }\quad \mathscr{D}=
\begin{pmatrix}0&1&D&D^2&\cdots&D^{n-1}\\
0&0&1&D&\cdots&D^{n-2}\\
0&0&0&1&\ddots&\vdots\\
\vdots&&&\ddots&\ddots&D\\
&&&&&1\\
0&&&\cdots&&0
\end{pmatrix}
\end{equation}
exactly as in (\ref{D-matrix}).  Moreover, 
\begin{equation}\label{condition-2}
v_{k-1}=0\quad\mbox{ if }\quad a^t\alpha^{k-1}=a^t\mathscr{D}\bigl(\varphi_k\alpha^k\bigr)=0. 
\end{equation}
And $v_n=0$ provided $\varphi_1\cdots\varphi_n=a_n$.  Thus we have shown that the construction of the auxiliary functions detailed in \S\ref{sec-auxiliary} ensures $v=0$.  It remains to verify that the initial conditions imposed on $\psi_{n,1}=\tau_n$ are satisfied. 

By (\ref{tau-derivatives}), 
\begin{equation}\label{tau-derivatives-at-0}
\begin{split}
\tau_n^{(j)}(0)&=(\tau_1(0),\ldots,\tau_n(0))\left(\left.(D+\Phi)^{j}\mathbf{e}^n\right|_{0}\right)\\
&=(\mathbf{e}^n)^t\left(\left.(D+\Phi)^{j}\mathbf{e}^n\right|_{0}\right),
\end{split}
\end{equation}
since (\ref{zero-values}) implies $(\tau_1(0),\ldots,\tau_n(0))=(\mathbf{e}^n)^t$. 
The quantity $(\mathbf{e}^n)^t\left(\left.(D+\Phi)^{j}\mathbf{e}^n\right|_{0}\right)$ can be read off the bottom row of $A$ has depicted in (\ref{A-form}), whereby 
\begin{equation}\nonumber
\tau_n^{(j)}(0)=A_{n\,n+1-j}=(\mathbf{e}^n)^t\left(\left.(D+\Phi)^{j}\mathbf{e}^n\right|_{0}\right)=0\quad\mbox{ if }\quad 1\leq j\leq n-1. 
\end{equation}
Since $\tau_n(0)=1$ by (\ref{zero-values}), this proves that $\psi_{n,1}=\tau_n$ satisfies the initial conditions stated in Theorem~\ref{thm-main}. 

\subsection{Proof of Theorem~\ref{thm-main}, part 2\label{proof-part-two}}

To begin we examine further the matrix $A$ defined by (\ref{A}), having the form (\ref{A-form}).  The goal is to determine $A(0)$ at every entry not in the first column. 
The auxiliary functions $\varphi_k$ are constructed according to initial conditions (\ref{k-initial}).  These initial conditions imply 
\begin{equation}\label{phi-initial}
\varphi_k(0)=1 \mbox{ for every }1\leq k\leq n\quad\mbox{ and }\quad \varphi_k^{(s)}(0)=0\mbox{ for every }1\leq s\leq k-2.
\end{equation}
It follows in turn from (\ref{A-form}) that $A_{j\,j+1}(0)=1$ for every $1\leq j\leq n$.  

We claim that for every $1\leq j\leq n$ and every $2\leq k\leq j-1$, $A_{jk}(0)=0$.  To see this, recall from \S\ref{sec-proof-part-one} that $A_{jk}$ is the sum of all terms of the form $(\mathbf{e}^j)^tw\mathbf{e}^n$, where 
\begin{equation}\nonumber 
w=w_1\cdots w_{n+1-k},\quad w_{n+1-k}=\Phi,\quad \rho_\Phi(w)=n-j\quad\mbox{ and }\quad \rho_D(w)=j+1-k, 
\end{equation}
the latter condition following from the fact that $\rho_\Phi(w)$ and $\rho_D(w)$ sum to the total degree $n+1-k$.  Since $2\leq k\leq j-1$ by assumption, it follows that 
\begin{equation}\label{rho-D}
1\leq \rho_D(w)\leq j-1.
\end{equation}
Thus each term of the form $(\mathbf{e}^j)^tw\mathbf{e}^n$ consists of $j+1-k$ occurrences of $D$, interspersed with $\varphi_{j+1},\ldots,\varphi_{n}$, the right-most term being $\varphi_n$.  Upon expansion by evaluation of the derivatives, such a term transforms to a sum of products
\begin{equation}\label{products}
\varphi_{j+1}^{(s_{j+1})}\cdots\varphi_n^{(s_n)}\quad\mbox{ where }\quad s_{j+1}+\cdots+s_n=\rho_D(w)=j+1-k. 
\end{equation}
Therefore at least one of the orders $s_r$ is positive, and every order satisfies 
\begin{equation}\nonumber
s_r\leq \rho_D(w)\leq j-1.
\end{equation}
Hence any positive order $s_r$ satisfies $1\leq s_r\leq r-2$. It follows from (\ref{phi-initial}) that $\varphi_{r}^{(s_r)}(0)=0$. The entry $A_{jk}$ being a sum of such products (\ref{products}), it follows that $A_{jk}(0)=0$ as claimed. 

To summarize, the matrix $A(0)$ has entries  
\begin{equation}\label{entries}
A_{j\,j+1}(0)=1\quad\mbox{ and }\quad A_{jk}(0)=0\mbox{ for }2\leq k\neq j+1.
\end{equation}

Returning to the statement of Theorem~\ref{thm-main}, let us now consider $\psi_{n,k}$ for some $k\neq 1$.  The analysis is very similar to that for $\psi_{n,1}$ carried out in \S\ref{sec-proof-part-one}, except that now we make use of the additional information concerning $A(0)$ obtained above.  Indeed, the analysis of $\psi_{n,1}$ did not involve the initial conditions imposed on the auxiliary functions. Those initial conditions now come into play. 

With $k$ fixed, set $\eta=\pi_n^{k-1}$ as in Theorem~\ref{thm-main}, and set 
\begin{equation}\nonumber
\sigma_n=\psi_{n,k}=T_{\varphi_{\eta(1)},\ldots,\varphi_{\eta(n)};\eta^{-1}(n)}. 
\end{equation}
For each $1\leq j\leq n-1$, set 
\begin{equation}\nonumber
\sigma_j=T_{\varphi_{\eta(1)},\ldots,\varphi_{\eta(n)};\eta^{-1}(j)}. 
\end{equation}
Observe by Lemma~\ref{lem-trig} that if 
\begin{equation}\nonumber
f=\sum_{j=1}^nc_j\sigma_j
\end{equation}
for some coefficient functions $c_j\in\acloc$, then 
\begin{equation}\label{iteration-2}
f^\prime=\sum_{j=1}^n\gamma_j\sigma_j\quad\mbox{ where }\quad \begin{pmatrix}\gamma_1\\ \vdots\\ \gamma_n\end{pmatrix}=(D+\Phi)\begin{pmatrix}c_1\\ \vdots\\ c_n\end{pmatrix},
\end{equation}
exactly as in (\ref{iteration}).  Hence precisely the same analysis yields
\begin{equation}\label{sigma-formula}
\sum_{j=0}^na_j\sigma_n^{(n-j)}=(\sigma_1,\ldots,\sigma_n)v
\end{equation}
where, as in \S\ref{sec-proof-part-one}, 
\begin{equation}\label{v}
v=\begin{pmatrix}v_1\\ \vdots\\ v_n\end{pmatrix}=\sum_{j=0}^na_j(D+\Phi)^{n-j}\mathbf{e}^n=Aa, 
\end{equation}
with $a$ and $A$ given by (\ref{a}) and (\ref{A}).  It follows that $v=0$ as before, and hence that $\psi_{n,k}=\sigma_n$ is a solution to (\ref{ode-n}).  

On the other hand, (\ref{zero-values}) implies 
\begin{equation}\nonumber
(\sigma_1(0),\ldots,\sigma_n(0))=\bigl(\mathbf{e}^{n-k+1}\bigr)^t.
\end{equation}
So, unlike (\ref{tau-derivatives-at-0}), for $0\leq j\leq n-1$, we have
\begin{equation}\label{sigma-derivatives-at-0}
\begin{split}
\sigma_n^{(j)}(0)&=(\sigma_1(0),\ldots,\sigma_n(0))\left(\left.(D+\Phi)^{j}\mathbf{e}^n\right|_{0}\right)\\
&=(\mathbf{e}^{n-k+1})^t\left(\left.(D+\Phi)^{j}\mathbf{e}^n\right|_{0}\right)\\
&=A_{n-k+1\,n+1-j}(0)\\
&=\left\{
\begin{array}{cc}
1&\mbox{ if } j=k-1\\
0&\mbox{ otherwise }
\end{array}
\right..
\end{split}
\end{equation}
Thus $\psi_{n,k}=\sigma_n$ satisfies precisely the initial conditions asserted in Theorem~\ref{thm-main}. This completes the proof.

\section{Remarks\label{sec-remarks}}

The general solution to (\ref{ode}) in the case $n=2$ first appeared in \cite{Gi:arXiv2025}; the result for higher order cases $n>2$ as described here in Theorem~\ref{thm-main} is new.  A qualitative difference between order $n=2$ and higher order $n>2$ is that the solution formula supplied by Theorem~\ref{thm-main} is global for second order equations, but in general only local for higher order.  Indeed, the hypothesis of Theorem~\ref{thm-main} that the coefficients belong to $\lloc$ could be relaxed in an obvious way.  Namely, one could assume just that the coefficients belong to $L^1(\Omega)$ for some open interval $\Omega\ni0$, and obtain precisely the same conclusion for equations of order $n>2$. No essential changes to the proofs given in \S\ref{sec-proof} are needed to accommodate this slight generalization.  The original hypothesis that coefficients belong to $\lloc$ has been retained for simplicity of exposition, and to accommodate the stronger conclusion valid for equations of order $n=2$. 

A second remark concerns the notions of explicit solution and solution by quadrature. Liouville famously proved in 1841 that Riccati equations of a specific type admit no solution by means of a certain limited class of operations (algebraic, exponential and logarithmic), with the consequence that associated second order linear ODEs also admit no such solution \cite{Li:1841}.  But this does not preclude solution by quadrature in the usual sense for general second order equations, since indefinite integration is not included among Liouville's limited class of operations, and Liouville himself 
points out 
that allowing integration would expand the class of solvable equations.  Thus Theorem~\ref{thm-main} does not run afoul of the classical impossibility result.  But does it constitute explicit solution by quadrature? Should the Dyson series solution (\ref{system-solution}) be considered explicit solution by quadrature?  This is ultimately a matter of convention, but it is worth reiterating that in the case $n=1$ the Dyson series reproduces the standard exponential solution $e^{\int_0^x\coeff{m}}$, which is universally accepted as explicit solution by quadrature of the corresponding first order equation.  Indeed, standard notation conceals the fact that the scalar exponential primitive is naturally expressed as a series.  Like the exponential primitive, trig operators involve series of integrals, and so it seems consistent to view them too as a sort of explicit quadrature, their novel formulation notwithstanding.  

In any case, Theorem~\ref{thm-main} provides a systematic, concrete formulation of the solution to an important class of equations, thereby filling a major gap in the elementary theory of differential equations. From an applied perspective, this offers the prospect of fresh insights in the various contexts in which linear ODEs arise. 

\appendix

\section{Convergence of the matrix Dyson series\label{sec-convergence}}

Fix an entrywise locally integrable matrix function $\mathbf{m}:\real\rightarrow\complex^{n\times n}$. We verify that the first order initial value problem in $n\times n$ matrices,  
\begin{equation}\label{system-2}
M^\prime=\coeff{m}M,\qquad M(0)=I
\end{equation}
has an explicit locally absolutely continuous solution for $x>0$ of the form
\begin{multline}\label{system-solution-2}
M(x)=I+\int\limits_{0<s_1<x}\coeff{m}(s_1)\,ds_1+\\
\int\limits_{0<s_1<s_2<x}\coeff{m}(s_2)\coeff{m}(s_1)\,ds_1ds_2+\int\limits_{0<s_1<s_2<s_3<x}\coeff{m}(s_3)\coeff{m}(s_2)\coeff{m}(s_1)\,ds_1ds_2ds_3+\cdots
\end{multline}

To begin, denote
\begin{equation}\label{denote}
\mathbf{M}_j(x)=\int\limits_{0<s_1<\cdots<s_{j}<x}\coeff{m}(s_{j})\coeff{m}(s_{j-1})\cdots\coeff{m}(s_1)\,ds_1\cdots ds_{j}
\end{equation}
and observe that 
\begin{equation}\label{re-write}
\mathbf{M}_{j+1}(x)=\int_0^x\coeff{m}(s)\mathbf{M}_j(s)\,ds.
\end{equation}
If $\mathbf{M}_j$ is locally absolutely continuous, then $\coeff{m}\mathbf{M}_j$ is locally integrable, and it follows by (\ref{re-write}) that $\mathbf{M}_{j+1}$ is locally absolutely continuous.  Thus by induction $\mathbf{M}_j$ is locally absolutely continuous on the positive half line for all $j\geq 1$, and 
\begin{equation}\label{term-differentiation}
\frac{d}{dx}\mathbf{M}_{j+1}(x)=\coeff{m}(x)\mathbf{M}_j(x).
\end{equation} 

For $1\leq k,\ell,\leq n$, let $\coeff{m}_{k\ell}$ denote the $(k,\ell)$-entry of $\coeff{m}$, and set 
\begin{equation}\label{g-max}
g=\max_{1\leq k,\ell\leq n}|\coeff{m}_{k\ell}|. 
\end{equation}
Then $g$ is locally integrable, and the modulus of each entry of the matrix
\[
\coeff{m}(s_j)\coeff{m}(s_{j-1})\cdots\coeff{m}(s_1)
\] 
is less than or equal to 
\begin{equation}\label{upper-bound}
n^{j-1}g(s_j)g(s_{j-1})\cdots g(s_1).
\end{equation}
Symmetry of the tensor product (\ref{upper-bound}) implies 
\begin{equation}\label{prod-g}
\int\limits_{0<s_1<\cdots<s_j<x}n^{j-1}g(s_j)g(s_{j-1})\cdots g(s_1)\,ds_1\cdots ds_j=\frac{1}{n}\frac{\left( n\int_0^xg\right)^j}{j!}.
\end{equation}
Therefore initial sums of the series 
\[
\sum_{j=1}^\infty\coeff{m}(x)\mathbf{M}_j(x)
\]
are entrywise uniformly dominated by 
\begin{equation}\label{dominator}
e^{n\int_0^xg}\coeff{m}(x). 
\end{equation}
Setting $\mathbf{M}_0=I$, it follows by dominated convergence that 
\begin{equation}\label{dominated-convergence}
\sum_{j=1}^\infty\mathbf{M}_j(x)=\sum_{j=0}^\infty\int_0^x\coeff{m}(s)\mathbf{M}_j(s)\,ds=\int_0^x\sum_{j=0}^\infty\coeff{m}(s)\mathbf{M}_j(s)\,ds
\end{equation}
Thus, referring to (\ref{system-solution-2}), $M(x)=I+\sum_{j=1}^\infty\mathbf{M}_j(x)$ is locally absolutely continuous, being the sum of $I$ and the primitive of an integrable function.  Term-by-term differentiation yields (\ref{system-2}).

\section{Parts I and II of Corollary~\ref{cor-examples} in detail\label{sec-derivation}}

As before, let $P:\lloc\rightarrow\acloc$ denote the primitive operator $Pf(x)=\int_{0}^xf$. 
\subsection{Second order equations\label{sec-second-order}}
Consider the equation 
\begin{equation}\label{second-order}
y^{\prime\prime}=a_1y^\prime+a_2y,
\end{equation}
where $a_1,a_2\in\lloc$ are given. 
The first step is to solve the auxiliary equations. Following \S\ref{sec-auxiliary}, $\alpha^3=(0,0,1)^t$.  To determine the auxiliary function $\varphi_2$ one solves 
\begin{equation}\label{first-aux-degree-2}
(-1,a_1,a_2)\begin{pmatrix}0&1&D\\ 0&0&1\\ 0&0&0\end{pmatrix}\bigl(u\alpha^3\bigr)=0\quad\Leftrightarrow\quad -u^\prime+a_1u=0
\end{equation}
with initial condition $u(0)=1$.  The solution is $u=\varphi_2=e^{Pa_1}$.  Therefore 
\begin{equation}\nonumber
\alpha^2=\begin{pmatrix}0&1&D\\ 0&0&1\\ 0&0&0\end{pmatrix}\bigl(\varphi_2\alpha^3\bigr)=\begin{pmatrix}a_1e^{Pa_1}\\ e^{Pa_1}\\ 0\end{pmatrix}.
\end{equation}
The remaining auxiliary function $u=\varphi_1$ by definition satisfies the equation
\begin{equation}\nonumber
(-1,a_1,a_2)\begin{pmatrix}0&1&D\\ 0&0&1\\ 0&0&0\end{pmatrix}\bigl(u\alpha^2\bigr)=-a_2\quad\Leftrightarrow\quad -ue^{Pa_1}=-a_2.
\end{equation}
Thus $\varphi_1=a_2e^{-Pa_1}$.  According to Theorem~\ref{thm-main}, the general solution to (\ref{second-order}) is
\begin{equation}\label{second-order-solution}
y=y(0)C+y^{\prime}(0)S\quad\mbox{ where }\quad C=T_{a_2e^{-Pa_1},e^{Pa_1};2}\mbox{ and }S=T_{e^{Pa_1},a_2e^{-Pa_1};1}.
\end{equation}
That (\ref{second-order-solution}) solves (\ref{second-order}) is easily verified directly using Lemma~\ref{lem-trig}.  The explicit solution (\ref{second-order-solution}) has numerous applications, leading to explicit integral formulas based on the multex series (\ref{multex}). 

\subsection{Third order equations\label{sec-third-order}}

The solution (\ref{second-order-solution}) to second order equations is needed to determine the auxiliary functions for the general third order equation
\begin{equation}\label{third-order}
y^{(3)}=a_1y^{\prime\prime}+a_2y^\prime+a_3y.
\end{equation}
Following the prescription in \S\ref{sec-auxiliary} as before, set $\alpha^4=(0,0,0,1)^t$, and let $u=\varphi_3$ denote the unique solution to 
\begin{equation}\nonumber
(-1,a_1,a_2,a_3)\begin{pmatrix}0&1&D&D^2\\ 0&0&1&D\\ 0&0&0&1\\ 0&0&0&0\end{pmatrix}\bigl(u\alpha^4\bigr)=0\quad\Leftrightarrow\quad -u^{\prime\prime}+a_1u^\prime+a_2u=0
\end{equation}
such that $u(0)=1$ and $u^\prime(0)=0$.  Thus $\varphi_3=C=T_{a_2e^{-Pa_1},e^{Pa_1};2}$ as in (\ref{second-order-solution}), and
\begin{equation}\nonumber
\alpha^3=\begin{pmatrix}0&1&D&D^2\\ 0&0&1&D\\ 0&0&0&1\\ 0&0&0&0\end{pmatrix}\bigl(\varphi_3\alpha^4\bigr)=
\begin{pmatrix}C^{\prime\prime}\\ C^\prime\\ C\\ 0\end{pmatrix}=
\begin{pmatrix}a_1e^{Pa_1}T_{a_2e^{-Pa_1},e^{Pa_1};1}+a_2C\\[5pt] e^{Pa_1}T_{a_2e^{-Pa_1},e^{Pa_1};1}\\[3pt] C\\ 0\end{pmatrix}.
\end{equation} 
The auxiliary function $u=\varphi_2$ is the solution to the equation
\begin{equation}\nonumber
(-1,a_1,a_2,a_3)\begin{pmatrix}0&1&D&D^2\\ 0&0&1&D\\ 0&0&0&1\\ 0&0&0&0\end{pmatrix}\bigl(u\alpha^3\bigr)=0\quad\Leftrightarrow\quad 
-uC^\prime-(uC)^\prime+a_1uC=0
\end{equation}
subject to the initial condition $u(0)=1$.  Thus
\begin{equation}\nonumber
\varphi_2=e^{P(a_1-2C^\prime/C)}=e^{Pa_1}C^{-2}\quad\mbox{ and }\quad 
\alpha^2=\begin{pmatrix}\varphi_2C^\prime+\bigl(\varphi_2C\bigr)^\prime\\ \varphi_2C\\ 0\\ 0\end{pmatrix}.
\end{equation}
The equation for $u=\varphi_1$ is then
\begin{equation}\nonumber
(-1,a_1,a_2,a_3)\begin{pmatrix}0&1&D&D^2\\ 0&0&1&D\\ 0&0&0&1\\ 0&0&0&0\end{pmatrix}\bigl(u\alpha^2\bigr)=-a_3
\quad\Leftrightarrow\quad-u\varphi_2C=-a_3.
\end{equation}
Since $C=\varphi_3$, the above equation implies $\varphi_1\varphi_2\varphi_3=a_3$, whereby $\varphi_1=a_3e^{-Pa_1}C$. 

Having determined $\varphi_1,\varphi_2,\varphi_3$, Theorem~\ref{thm-main} implies the general solution to (\ref{third-order}) is 
\begin{equation}\label{third-order-solution-2}
y=y(0)T_{\varphi_1,\varphi_2,\varphi_3;3}+y^\prime(0)T_{\varphi_3,\varphi_1,\varphi_2;1}+y^{\prime\prime}(0)T_{\varphi_2,\varphi_3,\varphi_1;2}
\end{equation}
where 
\begin{equation}\label{phi-3}
\varphi_1=a_3e^{-Pa_1}C,\quad\varphi_2=e^{Pa_1}C^{-2},\quad\mbox{ and }\quad\varphi_3=C=T_{a_2e^{-Pa_1},e^{Pa_1};2}.
\end{equation}
It is a straightforward if tedious exercise to check this directly using Lemma~\ref{lem-trig}.  

\section*{Acknowledgements}
Research was partially supported by NSERC grant DG RGPIN-2022-04547.

%


\begin{thebibliography}{10}

\bibitem{BeDeTo:1988}
R.~Beals, P.~Deift, and C.~Tomei.
\newblock {\em Direct and inverse scattering on the line}, volume~28 of {\em
  Mathematical Surveys and Monographs}.
\newblock American Mathematical Society, Providence, RI, 1988.

\bibitem{BoDi:1965}
W.~E. Boyce and R.~C. DiPrima.
\newblock {\em Elementary differential equations and boundary value problems}.
\newblock John Wiley \& Sons, Inc., New York-London-Sydney, 1965.

\bibitem{DoFr:1977}
J.~D. Dollard and C.~N. Friedman.
\newblock Product integrals and the {S}chr\"{o}dinger equation.
\newblock {\em J. Mathematical Phys.}, 18(8):1598--1607, 1977.

\bibitem{Dy:1949}
F.~J. Dyson.
\newblock The {S} matrix in quantum electrodynamics.
\newblock {\em Physical Review}, 75(11):1736--1755, Jun 1949.

\bibitem{encyclopedia_mathematics_linear_ode}
{Encyclopedia of Mathematics}.
\newblock Linear ordinary differential equation --- {Encyclopedia of
  Mathematics}.
\newblock
  \url{https://encyclopediaofmath.org/wiki/Linear_ordinary_differential_equation},
  2020.
\newblock Last updated: 2020-06-05. Accessed: 2025-12-07.

\bibitem{Gi:SIMA2024}
P.~C. Gibson.
\newblock Explicit solution of the 1{D} {S}chr\"odinger equation.
\newblock {\em SIAM J. Math. Anal.}, 56(4):4466--4493, 2024.

\bibitem{Gi:arXiv2025}
P.~C. Gibson.
\newblock Solution of the scalar {R}iccati equation, 2025.
\newblock arXiv preprint 2508.02784, 22 pages.

\bibitem{Gi:Proc2026}
P.~C. Gibson.
\newblock On integration of the scalar {R}iccati equation.
\newblock {\em Proceedings of the American Mathematical Society, Series A}, in
  press.
\newblock Published electronically.

\bibitem{GiJo:1990}
R.~D. Gill and S.~Johansen.
\newblock A survey of product-integration with a view toward application in
  survival analysis.
\newblock {\em Ann. Statist.}, 18(4):1501--1555, 1990.

\bibitem{GrNg:2019}
E.~Grenier and T.~T. Nguyen.
\newblock Green function of {O}rr-{S}ommerfeld equations away from critical
  layers.
\newblock {\em SIAM J. Math. Anal.}, 51(2):1279--1296, 2019.

\bibitem{He:1987}
I.~H. Herron.
\newblock The {O}rr-{S}ommerfeld equation on infinite intervals.
\newblock {\em SIAM Rev.}, 29(4):597--620, 1987.

\bibitem{In:1944}
E.~L. Ince.
\newblock {\em Ordinary {D}ifferential {E}quations}.
\newblock Dover Publications, New York, 1944.

\bibitem{La:1869}
J.-L. Lagrange.
\newblock Sur les int\'egrales particuli\`eres des \'equations
  diff\'erentielles.
\newblock In J.-A. Serret, editor, {\em \OE uvres de Lagrange}, volume~4, pages
  5--108. Gauthier-Villars, Paris, 1869.
\newblock Reprint of the 1774 memoir.

\bibitem{LaTr:1997}
H.~Langer and C.~Tretter.
\newblock Spectral properties of the {O}rr-{S}ommerfeld problem.
\newblock {\em Proc. Roy. Soc. Edinburgh Sect. A}, 127(6):1245--1261, 1997.

\bibitem{Li:1841}
J.~Liouville.
\newblock Remarques nouvelles sur l’équation de {R}iccati.
\newblock {\em J. Math. Pures Appl., Ser. 1}, 6, 1841.

\bibitem{Pe:1888}
G.~Peano.
\newblock Int\'{e}gration par s\'{e}ries des \'{e}quations diff\'{e}rentielles
  lin\'{e}aires.
\newblock {\em Math. Ann.}, 32(3):450--456, 1888.

\bibitem{PiMo:2025}
A.~Piterskaya and M.~Mortensen.
\newblock A study of the {O}rr-{S}ommerfeld and induction equations by
  {G}alerkin and {P}etrov-{G}alerkin spectral methods utilizing {C}hebyshev
  polynomials.
\newblock {\em J. Comput. Appl. Math.}, 459:Paper No. 116374, 14, 2025.

\bibitem{SaNa:2021}
J.~J. Sakurai and J.~Napolitano.
\newblock {\em Modern quantum mechanics}.
\newblock Cambridge University Press, Cambridge, United Kingdom ;, third
  edition. edition, 2021 - 2021.

\bibitem{Sl:2007}
A.~Slav\'{\i}k.
\newblock {\em Product integration, its history and applications}, volume~1 of
  {\em Ne\u{c}as Center for Mathematical Modeling}.
\newblock Matfyzpress, Prague, 2007.
\newblock Corrected translation of the Czech original, D\v{e}jiny
  Matematiky/History of Mathematics, 29.

\bibitem{Te:2012}
G.~Teschl.
\newblock {\em Ordinary differential equations and dynamical systems}, volume
  140 of {\em Graduate Studies in Mathematics}.
\newblock American Mathematical Society, Providence, RI, 2012.

\bibitem{HoVo:1938}
V.~Volterra and B.~Hostinsky.
\newblock {\em Operations infinitesimales lineaires.}
\newblock Collection de monographies sur la theorie des fonctions.
  Gauthier-Villars, Paris, 1938.

\end{thebibliography}

\end{document}